\documentclass{paper}[12pt] %For draft work

\usepackage{amssymb}   
\usepackage{amsthm}
\usepackage{amsmath}
\usepackage{amsfonts}
\usepackage{latexsym} 
\usepackage{eucal}
\usepackage{indentfirst} 
\usepackage{graphicx} %Allow use of pictures
\usepackage{verbatim} %Allow \begin{comment} and \end{comment} 

\def\R{{\mathbb R}} % Reals
\begin{document}

%\begin{document}
\overfullrule=0pt
\baselineskip=24pt
% if i uncomment the preceeding line, the text will be double-spaced.
\font\tfont= cmbx10 scaled \magstep3
\font\sfont= cmbx10 scaled \magstep2
\font\afont= cmcsc10 scaled \magstep2
\title{ }
\bigskip
\bigskip
\author{{\it Differential Forms: Theory and Practice.} \\
By Steven H. Weintraub, Academic Press, 2014, \\xi +  395 pp., ISBN 9780123944030, \$119.95\\ \\
{\it Reviewed by} Thomas Garrity}

\date{}

\maketitle

Though differential forms are an essential tool for much of modern mathematics, they are not standardly taught at the undergraduate level.  In fact, they are 
far from being well-known by many working mathematicians, as I discovered while writing this review, when     I sent out the following to a number of mathematician friends:

``What I am wondering is how much are differential forms  common knowledge among professional mathematicians.  I seem to recall that differential forms were assumed knowledge in my first year in graduate school (an incorrect assumption for most if not all of the incoming graduate students).  For some of you, differential forms are tools you use almost daily.  That is the case for me.  For some, differential forms are something you sort of remember from graduate school. And possibly for some, they are mere words, if that.''

The replies ranged from, ``I use them everyday,'' to `` occasionally they come up in a paper I am reading or writing, in which case I need to relearn the basics,'' to ``dimly recalled from a physics course from years ago,'' to ``mere words.''

\section{The Pedagogical Challenge }
Among the basic questions that need to be asked to understand the world are 
\begin{enumerate}
\item How much is there?
\item How fast is it changing?
\end{enumerate}
To some extent single variable and multi-variable calculus give us the tools to answer these questions, with the integral answering a ``how much is there'' question and the derivative answering a ``how fast is it changing'' question.  The close to miraculous Fundamental Theorem of Calculus asserts that the answers to these two seemingly distinct questions are closely and intimately linked. 

Differential forms can be interpreted as providing a means for answering these questions in more general contexts.  And as with calculus, there is a corresponding fundamental theorem called Stokes' Theorem. The differential forms version of Stoke's theorem has as special cases the Divergence Theorem, Green's Theorem and Stokes' Theorem for surfaces.

The pedagogical challenges in teaching differential forms parallel the teaching challenges for calculus. First off, both the derivative and the integral have meaning.  In a beginning course, usually these meanings are sketched out at an intuitive level.  In most beginning real analysis courses, the more rigorous foundations behind the meanings are laid out. But a large part of a beginning course concentrates not on the meaning but on computing derivatives and integrals.  The reason calculus is so successful is  both that its basic terms have meaning and that they can be easily computed.

But knowing the meaning of the derivative  is quite different from the actual practice of computing the derivative.  After all, knowing the derivative is a measure of a rate of change of a function helps little in computing the derivative of $ {\it cos\left(x^3 + 4x^2 + 5x + 9 \right)}$.  All of us would simply use our knowledge of the chain rule, of the derivative of cosine and of the  derivative of a polynomial.  In fact, consciously thinking of the meaning would get in the way of getting the answer.  The computational side is quite distinct from the meaning side.

The same difficulty holds in the teaching of  differential forms.   Their meaning is tied to measuring.  To some extent, a differential $k$-form $\omega$ on $\R^n$ is a device that will allow the measuring of the volume of a dimension $k$ submanifold $M$.  The technical definitions are to some extent getting to the point where the symbol
$$\int_M \omega$$
makes sense as a number that has properties that one would expect of a volume measurement.  Similarly, their meaning is also tied to rates of change.  Here we need to carefully define what taking a derivative of a $k$-form could possibly mean.  Thus we need to make sense out of the symbol
$$\mathrm{d} \omega.$$
The operation of going from a $k$-form $\omega$  to the $(k+1)$-form $\mathrm{d}\omega$ is called taking the exterior derivative. 
 The link with the measuring of volume is the true punchline in the beginning theory of differential forms, namely Stokes' Theorem:
$$\int_{\partial (M)} \omega    = \int_M \mathrm{d} \omega .$$
Here we are letting $M$ be a $(k+1)$-dimensional manifold, with the new symbol $\partial(M)$ denoting the boundary of $M$.  Thus $\partial(M)$ is itself $k$-dimensional.  A more intuitive interpretation is to think of the initial $\omega$ as analogous to a function (which, to be clear, it is not in general) and the taking of the exterior derivative as the taking of a derivative.  Then Stokes' Theorem can be described as 
$$\mbox{ The average of a function over the boundary of a manifold}$$
equals 
$$\mbox{ The average of the derivative over a manifold}.$$
Making this rigorous takes work (and should be viewed as a true triumph of early twentieth century mathematics \cite{Hawkins05,   Cogliati11}), involving not only definitions of differential forms and exterior derivatives but also the definitions of manifolds and boundaries.  Once the correct definitions are given, though, the proof of Stokes' Theorem should become straightfoward.

But none of these definitions particularly help in calculating.  Here one must simply practice, practice, practice, just as most of us once did with multiplication tables and the taking of one-variable derivatives.

I suspect that one reason some of the mathematicians whom I asked about differential forms had only  a hazy understanding, at best, was that they did not practice when they should have.  Probably they were quickly introduced to differential forms in graduate school.  At this level of mathematical maturity, the rigorous development can be done fairly quickly, after which the professor would say something along the lines of ``The rest is just calculations.''
This is a true statement, save for the use of the word ``just.''  Many graduate students hear that word and then don't  bother to ``practice their scales'' and do the computations.   (Most of us have seen first year college students in beginning calculus who say that they know what  a derivative is but who can't compute any derivative beyond $x^2$;  this is the graduate student analog.)

Thus the teaching of differential forms must explicitly discuss their underlying meanings while at the same time must emphasize the computational aspects.  

\section{Why Teach Differential Forms, or, Why They are Important}

There are two questions here, namely, why teach differential forms at all and why teach them to undergraduates.  The why for teaching differential forms at all is the same as asking why differential forms are important in mathematics. Here the answer comes down to the fact that differential forms are to a large extent the best tools for answering questions about how much is there (questions of area and volume) and questions about rate of change. Considering how much of mathematics can be phrased as answering these two questions, the importance of differential forms is self-evident.  Let us look though at three examples, from seemingly three quite different branches of mathematics.  (There are many more.)

Differential geometry is the study of how geometric objects bend and twist. Slightly more technically, differential geometry can in part be taken to be the study of curvature in its many manifestations.  But curvature, to some extent, is measuring the rate of change of directions of tangents.  Once you invoke the phrase ``rate of change,'' differential forms must be relevant.  And for curvature, they are vital.

Second, consider  the four vector-valued partial differential equations that form the  Maxwell Equations.   These equations explicitly describe the relation between the electric field, the magnetic field, current and charge. Further, these equations lead to the recognition that light is a type of electromagnetic wave. In fact, from these equations you can deduce that  the speed of these electromagnetic waves is weird (the speed is the same for  observers in different frames of reference, as long as the observers are moving at constant velocity with respect to each other), which in turn leads directly to the special theory of relativity.  (For more on this, see the  recent text \cite{Garrity15}.)  But the electric field and the magnetic field can be encoded into a single differential two-form, called the electromagnetic two-form.  Maxwell's equations can then be described by how the exterior derivative acts on the electromagnetic two-form (and the Hodge star of the electromagnetic two-form).  Though initially this looks like a simple repackaging of the four original partial differential equations into the language of differential forms, this language allows for profound generalizations.  In part, it will eventually lead to recognizing that the forces coming from the electric and magnetic fields can actually be described as a type of curvature, leading in turn to one of the themes of 20th century physics:
$$\mbox{Force}=\mbox{Curvature}.$$

Finally let us turn to numerical analysis.  Certainly a large part of numerical analysis is concerned with approximating areas and volumes and approximating rates of change.  I suspect, though, that most traditional numerical analysts would not bother with the seemingly abstract nature of differential forms. This view is currently changing.  See Douglas Arnold's Plenary Address at the 2002 International Congress of Mathematics \cite{Arnold02} , his  important, book length survey paper \cite{Arnold-Falk-Winther06} (with coauthors Falk and Winther), and with the same coauthors, his \cite{Arnold-Falk-Winther10}.

Thus a professional mathematician should know about differential forms.  But how important are they for undergraduates to learn?   First for a caveat.  Given the choice, people should teach what they are mathematically excited about.  This is what fosters  for great teaching.   Many of us, though, find most of mathematics potentially exciting. Thus the choice of what we teach (again, in ideal world where we have total freedom in choosing what we teach) should be decided by what will help our students the most.  Most of our students will not become professional mathematicians.  For them, the most important thing is to take hard challenging courses that push them to a new level of mathematical maturity.  To some extent, the actual content is not that important.  Thus it is perfectly reasonable, if not desirable, to choose what to teach based on what will benefit those few students who will become professional mathematicians.  Based on this line of thought, the teaching of differential forms becomes eminently reasonable.  Another factor in determining what will help the future mathematician is not just what will be later important to know but also what is harder to ``learn in the streets.''  Most of the mathematics any practicing mathematician knows is, of course, not learned in a class.  Some topics are intrinsically easier to learn on their own, while for others, it is far better to learn in a class, with the rhythm of learning set by the class and instructor.  This is why it makes strong pedagogical sense to require   our majors to take a beginning real analysis course and a beginning course in group theory.  It would also justify the teaching of beginning representation theory to our majors, which is not that often done.  And it also justifies the teaching of differential forms, and in fact the teaching of the underlying meaning via rigorous definitions and by the requiring of many concrete calculational problems.  This is the goal of Weintraub's test.

\section{Weintraub's Approach}

In 1996, Steven Weintraub wrote {\it Differenential Forms:  A Complement to Vector Calculus} \cite{Weintraub96} with the aim of  introducing differential forms to students who have just finished a standard third-semester calculus course (vector calculus).  His current book, with its similar title, while built out of that earlier text, is aimed at more advanced students, namely students who have additionally taken linear algebra and who know the basics of point-set topology. In fact, I would be willing, if not eager, to use this as a text for a course with only a linear algebra prerequisite, filling in the background of the little point-set topology that is needed (mostly about open balls in $\R^n$).  

The first two chapters develop the basics of differential forms in $\R^n.$   There is a clear emphasis, especially in the exercises, on getting students to know how to manipulate and to  calculate with differential forms.  Chapter 3 moves up the level of abstraction a bit, though still with a solid computational feel.  The goal of this chapter is to show, for a differentiable map $f:\R^n \rightarrow \R^m,$ that $f$ sends points and tangent vector in $\R^n$ to points and tangent vectors in $\R^m,$ respectively, and show there is natural map $f^*$ that will send functions and differential forms on the initial range $\R^m$ to functions and differential forms on the original domain $\R^n,$ respectively. 

Chapter 4 is a  development of manifold theory.  Weintraub is careful with the needed subtleties, such as questions about orientation, boundaries and homotopies.
Building on Chapter 4's work on manifold theory, Chapter 5 develops vector bundles (in particular,  tangent and cotangent bundles) on manifolds. 
This allows him in Chapter 6  to define what it means to integrate a $k$-form along a $k$-dimensional manfold.

The last two chapters are his two main applications, or perhaps more accurately, his two main punchlines. Chapter 7 is a full proof of Stokes Theorem (a truly remarkable theorem) while Chapter 8 develops de Rham cohomology.

\section{Other Texts}
As differential forms are important, there are of course a number of texts in which they are introduced, including many beginning books on differential geometry and general relativity.  Here we will consider  a few texts which directly introduce differential geometry to undergraduates.

In 1968, Michael Spivak wrote his {\it Calculus on Manifolds} \cite{Spivak71}.  This classic text is a succinct and fairly abstract approach to differential forms  and is an excellent source for the mathematically mature.  He is quite explicit in following the philosophy that the bulk of the work should be in setting up the correct definitions.  Probably, though, it is not the right place for undergraduates to first learn about differential forms.  

From 1977, there is 
Schreiber's   {\it Differential Forms: A Heuristic Introduction} \cite{Schreiber84}, which is aimed at a lower level of mathematical maturity than Spivak or Weintraub, namely, for students who  know just vector calculus and a bit about  $\R^n$.   Schreiber is quite explicit that the exposition aims at the concrete.  For example, all of the work is done on submanifolds of $\R^n$.  One disadvantage is that there are no exercises.

My personal favorite  is   {\it Vector Calculus, Linear Algebra, and Differential Forms: A Unified Approach} \cite{Hubbard-Hubbard01}  by Hubbard and Hubbard.   To cover this book would require two semesters, which at most schools is not possible, which is the main reason I have never used this as a text for either multivariable calculus or for linear algebra.  This is unfortunate, as it does a wonderful job covering both meaning and computation.

R. Darling's {\it Differential Forms and Connections} \cite{Darling94} is aimed at a higher level of mathematical maturity, requiring a background in both mulitvariable calculus and linear algebra.  It moves at a quicker pace than Weintraub and as a consequence covers more material.  This is in part reflected in the title containing the term ``connection.''  This speed allows Darling to get to the  Yang-Mills equations and gauge theory by the end.

De Carmo's  {\it Differential Forms and Applications}
\cite{doCarmo00} from 2000 is a good book for senior math majors.  It is short, and quickly gets to the point of differential forms.

The  2012 edition of David Bachman's  {\it A Geometric Approach to Differential Forms }  \cite{Bachman11} is also quite good, and is comparable to Weintraub's text.  It is pitched to almost exactly the same audience.

The recent   {\it Exterior Analysis: Using Applications of Differential Forms } by \c{S}uhubi \cite{Suhubi13} is aimed at students at a slightly more advanced level of mathematical maturity.  It also frames the development in the language of the exterior algebra, which gives the book more of a linear algebra feel.  

Sinha's {\it Calculus of Tensors and Differential Forms} \cite{Sinha14} emphasizes the calculational side of differential forms, putting them into the slightly more general language of tensors.  The text is full of a massive number of calculations, with a plethora of indices.  If you want to emphasize the computations, this would not be a bad text.  Unfortunately, it has no exercises. To some extent, its approach is perpendicular to that of Spivak.

\section{Conclusion}
I liked Weintraub's book and look forward to teaching from it someday.   Also, if you are one of the mathematicians who have  only a hazy memory of differential forms, if that, then this would be an good book to go through on your own at night, working the exercises.  You will have an enjoyable time.

\noindent \it {Williams College, Williamstown MA 01267}

\noindent \it{tgarrity@williams.edu}

\end{document}